\documentclass[11pt]{article}
\usepackage{xcolor}
\usepackage{amsmath}
\usepackage{amssymb}
\usepackage{graphicx}
\usepackage{soul}
\usepackage{pgf}
\usepackage{tikz}
\usepackage{mathrsfs}

\newtheorem{theorem}{Theorem}[section]
\newtheorem{lemma}{Lemma}[theorem]
\newtheorem{proposition}{Proposition}

\newtheorem{corl}{Corollary}[lemma]
\newtheorem{rem}{Remark}[section]

\newtheorem{remth}{Remark}[theorem]
\newtheorem{note}{Note}[section]

\newcommand{\ix}{\underline{x}}          
\newcommand{\wt}[1]{\widetilde{#1}}
\newcommand{\xx}{[\![x]\!]}              

\newcommand{\SG}{{S}}                    
\newcommand{\SX}{{X}}                    
\newcommand{\SY}{{Y}}                    
\newcommand{\ovS}{\overline{\SG}}        
\newcommand{\symS}{\wt{\SG}}             
\newcommand{\xth}{\kappa}                   

\newcommand{\Frac}{\displaystyle\frac}            
\newcommand{\Int}{\displaystyle\int}              
\newcommand{\Max}{\displaystyle\max}              
\newcommand{\Sup}{\displaystyle\sup}              
\newcommand{\Sum}{\displaystyle\sum}              
\newcommand{\refrm}[1]{{\rm(\ref{#1})}}           

\newcommand{\vspandex}{\rule[-18pt]{0pt}{36pt}}   
\newcommand{\vspandexsmall}{\rule[-11pt]{0pt}{22pt}}   

\newcommand{\defif}[4]{
\left\{\begin{array}{ll}
#1,&\mbox{if $#2$},\\
#3,&\mbox{if $#4$}.\\
\end{array}\right.}

\newcommand{\defifsmall}[4]{
\left\{\begin{array}{ll}
\mbox{\small $#1$},&\mbox{\small if $#2$},\\
\mbox{\small $#3$},&\mbox{\small if $#4$}.\\
\end{array}\right.}

\newcommand{\be}{\begin{equation}}
\newcommand{\ee}{\end{equation}}
\newcommand{\Proof}{\noindent {\bf Proof.~}}
\newcommand{\QED}{\hfill \rule{2.5mm}{2.5mm}              

\vspace{7pt}                                              

\noindent }
\newcommand{\set}[1]{\left\{#1\right\} }          
\newcommand{\ex}[1]{\,{\rm exp}\set{#1}\,}        
\newcommand{\df}{:=}
\newcommand{\D}{\,\stackrel{\mbox{\footnotesize {\cal D}}}{=}\,}
\newcommand{\tss}[1]{\textsuperscript{#1}}       
\newcommand{\vv}{\vspace{7pt}

}
\DeclareMathOperator {\sign}{sign}

\newcommand{\signof}[1]{\sign(#1)}
\DeclareMathOperator{\supp}{supp}

\newcommand{\nn}{n}    

\def\ca#1{{\cal #1}}
\newcommand{\E}{{\mathbb E}}               
\newcommand{\PR}{{\mathbb P}}              

\newcommand{\BK}{{\mathbb K}}
\newcommand{\BL}{{\mathbb L}}
\newcommand{\N}{{\mathbb N}}
\newcommand{\R}{{\mathbb R}}
\newcommand{\BT}{{\mathbb T}}
\newcommand{\Var}{\text{\sf Var}}

\newcommand{\cS}{\mathscr{S}}              
\newcommand{\cM}{\mathscr{M}}              
\newcommand{\wek}[1]{\textsf{$\textbf{#1}$}}      

\newcommand{\I}[1]{1\hspace{-3pt} {\rm I}_{#1}}       
\newcommand{\Ib}[1]{1\hspace{-3pt} {\rm I}_{\{#1\}}}  

\begin{document}\sf
\hfill \today
\begin{center}
{\large\bf Random Gamma time}
\vv
Jerzy Szulga

Department of Mathematics and Statistics, Auburn University, USA
\vv\vv

\begin{minipage}{0.9\linewidth}
{\small {\bf Abstract.} We  discuss the  Gamma L\'evy process, including path properties, the inverse process, integrability, and its spin-offs obtained by compounding, exponentiation, and other operations; further extendable to arbitrary $\sigma$-finite continuous Borel spaces. An appendix on modular spaces and deterministic jump processes is included.}
\end{minipage}
\end{center}

{\footnotesize\tableofcontents}
\vv
\vv

\addcontentsline{toc}{section}{Preamble}
{\bf \Large Preamble}
\vv
A stochastic process $(\SX_t; t\in \R)$ entails algebraically the obvious integral $\SX f\df \int f\,d\SX$ for interval-based simple functions;  first with deterministic,  then with random values.  Extensions require an adequate structure of the process.
Origins of such approach can be traced back to \cite{Wie} and its temporaries.
\vv
A thorough description of arising ``integrable functions'' is one of the primary tasks. E.g., cf.\ \cite{KwaW},  for a L\'evy process, a derived random measure $\SX \I A$ would allow simple functions beyond just intervals, then  deterministic integrands that form a modular F-space (metric, complete, translation invariant). The three components (deterministic, diffusive, and Poissonian) jointly entail a rather sophisticated metric (cf.\ \cite[Th.\ 8.3.1]{KwaW}).
\vv
The diffusive part essentially falls into the Hilbertian theory but the diffusion-free part alone has its own specific impact often clouded and dominated by the present diffusion. In contrast, the study \cite{KalS} of pure Poissonian integrals  was confined to positive, centered, or symmetric pure jump processes. A domain $T\subseteq\R_+$ with the Lebesgue measure $\lambda$ suffices to derive an exact analysis of existence, convergence, or divergence; a blueprint ready to be transferable to more general Borel spaces $(\BL,\ca L,\lambda)$ with atomless $\sigma$-finite Borel measures $\lambda$ such as the Lebesgue measure on $\R^d$.
\vv
In contrast to general treatises on L\'evy process such as \cite{Ros,KalS,App}, etc.,  where particular processes serve merely as illustrations of the theory, we study the Gamma processes in detail, sort of like ``laboratory white mice''. Also for that reason we use the unit intensity $\lambda=1$ to avoid an unnecessary clutter; a non-unit intensity can be reinserted back if so desired. 
 Yet, we mark features and  characteristics that can be either  extended to general L\'evy process or are applicable to other specific cases such as stable processes. However, the Gamma process admits characteristics that are absent in many L\'evy process, e.g., moments of arbitrary order which, in turn, prompts for a study of their asymptotic behavior.
\vv
The first section recalls the Gamma process that generalizes arrival moments ${\SG}_n$ in a Poisson process  interpreted as a process of fractional signals on $[0,\infty)$,  expanding  the blueprint from \cite{KalS}. We note a connection to Thorin's GGC (Generalized Gamma Convolutions) (cf.\ \cite{Tho,Bar-NMS, BehB,PerS,JamRY} and numerous references therein);  however, positive Gamma integrals were considered there as individual entities rather than values of a function-indexed stochastic process. Gamma process may be composed at jumps with independent ``rewards'', including Bernoulli random variables that in the symmetric case entails a symmetric Gamma process.
\vv
The second section describes paths and $p$-moments and adds hyper-exponential moments of the inverse Gamma process; 
the appendix contains a primer of modular space augmented by basic deterministic calculus of jumps.
\vv
Following the typography of \cite{KalS}, for a general domain  $(\BL,\ca L,\lambda)$ we utilize the ``operator notation'' for integrals, e.g.,   $\lambda f\I T= \int_T f\,d\lambda$ or $\SG_T f=\int_T f\,d\SG$,   $T\in\ca L$. We may skip the subscript when the set is the entire space.
\vv
 When $\BL=[0,\infty)$  we specify $\lambda$ as the Lebesgue measure and with the linear order we further abbreviate  $\lambda_t f=\int_0^t f(x)\,dx$  and similarly ${\SG}_tf=\int_0^t f\,d{\SG}$, etc.
\section{The Gamma integral}
\subsection{Distributional approach}
The standard procedure leading to the so called ``stochastic integral'' is quite simple due to the specificity of the Gamma distribution. That is, the densities and Laplace transforms
\[
f_t(x)=\frac{x^{t-1}} {\Gamma(t) }\,e^{-x},\,x>0,\quad L_t(\theta)=\frac{1}{(1+\theta)^t}=e^{-t\ln(1+\theta)},\quad t,\theta >0
\]
yield finite dimensional distributions which by the Kolmogorov Consistency Theorem ensure the existence of a L\'evy process ${\SG}_t$  on $(0,\infty)$ with independent and stationary increments, augmented by ${\SG}_0=0$. For a real function $f\ge 0$ on $\R_+$, the following  function $L(\theta)$ is a Laplace transform of some probability measure that may be represented by a single random variable denoted, say, by $\SG _f$:
\be\label{Lf}
L(\theta)\df e^{-\lambda \ln (1+\theta f)} =:\E e^{-\theta \SG_f}.
\ee
The passage from an``entity'', i,e., a single random variable $\SG_f$ marked by $f$, to the linear process $\SG f$ requires some additional arguments. Then the variable $\theta$ becomes obsolete because it is built in the integrand.
Alternatively, we may utilize the Fourier transform
\be\label{Ff}
F(\theta)\df\E e^{\imath \theta \SG_f}=
\ex{-\frac 1 2 \lambda (1+\theta^2f^2)+\imath \,\lambda \arctan(\theta f)},
\ee
which stems from the L\'evy-Khinchin representation.
\vv
The probability distribution of $\SG _f$ is well defined on an arbitrary $(\BL, \ca L, \lambda)$, not just on the real half-line. Hence, we may consider the product space $(\BL',\ca L',\lambda') \df (\BK,\ca K,\kappa)\otimes (\BL,\ca L,\lambda)$. For example, consider a  probability distribution  $\kappa$ of a random variable $K$. The forthcoming Fourier transform describes first a random entity $\SG^{(\kappa)}_f$, then the linear stochastic ``compound'' or ``reward'' process $\SG^{(\kappa)} f$ (see also \refrm{SSf} below). Factually, there is no need to restrict $\kappa$ to be only a probability measure; it could be any measure satisfying the underlying integrability conditions.
\vv
In the simplest case, for a $(\pm 1)$-Bernoulli $K^{(\beta)}$, i.e.,  $\PR(K=1)=\alpha\in [0,1]$ and $\beta\df \E K=2\alpha-1\in [-1,1]$, we obtain a random  entity $\SG_f^{(\beta)}$:
\be\label{FfK}
\begin{array}{rclcl}
F^{(\kappa)}(\theta)&\df&\E e^{\imath \theta\SG^{(\kappa)}_f}&=&
\ex{-\frac 1 2 \lambda\, \E\, (1+\theta^2f^2)+
\imath \,\lambda\,\E\, \arctan(\theta K f)}\\
F^{(\beta)}(\theta)&\df&\E e^{\imath \theta \SG^{(\beta)}_f}&=&
\ex{-\frac 1 2 \lambda (1+\theta^2f^2)+\imath \,\beta \lambda \arctan(\theta f)}\\
\end{array}
\ee
where the upper mark ``$(\beta)$'' represents just the parameter of the Bernoulli distribution; e.g., in the symmetric case reduced to
$\symS_t\df \SG^{(0)}$ which also may be represented as $ {\SG}_{t/2}-{\SG}'_{t/2}$ with an independent subtrahend. 
\vv
In such approach the stochastic process $\SG_t$ is first acquired  only on $\R_+$. Then indeed it yields a linear stochastic process $\SG f\df \SG_f$, well defined on the set $\cS$ of step functions; and consequently,
\[
\SG f =\int_0^\infty f\,d\SG_t,\quad f\in \cS.
\]
Again, there is no need to specify the scalar $\theta$.
The range  $\SG (\cS)$ is a vector subspace of $L^0(\Omega)$. In view of the exponent, the set of step functions is dense in the modular metric space $L^{\phi_1}(\R_+)$ (cf.\ App.\ \ref{modspaces}).
\vv
It remains to see that the closure in $L^0$ (i.e., in probability) of  the core range corresponds to the whole $L^{\phi_1}$,  or to $L^{\phi_2}$ in the latter case. We will use the following Chebyshev's style inequality, involving just a single random variable $\SG_f$ rather than the process $\SG f$.
\begin{proposition} For a positive $\epsilon\le \frac 1 2$ and $f\in \cS$,
\[
 \PR\left(|\SG^{(\beta)} _f|>\epsilon\right)\le \left\{
\begin{array}{ll}
\vspandexsmall
\frac 3{\epsilon}\,  \phi_1(f), &\text{if } \beta\neq 0;\\
\frac 1{\epsilon}\,\phi_2(f),&\text{if }\beta=0.
\end{array}\right.
\]
\end{proposition}
\Proof
Put $X=\SG^{(\beta)}_f$. We apply the estimate (cf.\ \cite[Lemma 5.1(1)]{Kal},
\[
c(\epsilon)\df\PR(|X|>\epsilon)\le \Frac{\epsilon}{2} \Int_{-1/\epsilon}^{1/\epsilon}
\left|1-\E e^{\imath \theta X}\right|\,d\theta.
\]
Writing  $A+\imath \,B$ for the exponent  and
computing the modulus in the integrand
with $x=\cos B$ and $y=e^{-A}$,
\[
\left(1-2xy+y^2\right)^{1/2} \le \left((1-y)^2+2(1-x)\right)^{1/2}\le A+|B|.
\]
Both even functions, increasing on $\R_+$, are bounded by their end values. With $r=1/\epsilon$, we consider
\[
\Sup_x \Frac{\phi'_2(rx)}{\phi'_1(x)} \quad\text{and}\quad
\Sup_x \Frac{\phi'_0(rx)}{\phi'_1(x)}.
\]
So,
\[
\begin{array}{rclcl}
\vspandexsmall
\Frac{r^2(x+x^2)}{1+r^2x^2}&=&
r\,\Frac{r x}{1+r^2x^2}+\Frac{(rx)^2}{1+r^2x^2}
&\le& \frac r 2 +1\le r;\\
\Frac{r(1+x)}{1+r^2x^2}&=&
r\,\Frac{1}{1+r^2x^2}+\Frac{rx}{1+r^2x^2}&\le&  r +\frac 1 2\le \frac {3r} 2\\
\end{array}
\]
and adding when $\beta\neq 0$: $r+\frac{3r} 2\le 3r$.
\QED

Let us remind that thus far $\SG f$ has been well defined only for step functions on $\R_+$.
\begin{theorem}
On $\R_+$, the linear stochastic integral $\SG^{(\beta)} f$ is well defined on $L{^\phi_1}$ when $\beta\neq 0 $ and on $L^{\phi_2}$ when $\beta=0$.
\end{theorem}
\Proof
In both cases the closure argument is identical, so consider either $L^{\phi}$. For $f\in L^{\phi}$ there exists a sequence $f_n$ of step functions such that $f_n\to f$ in $L^\phi$, i.e., $f_n$ is Cauchy. By the above Proposition,  $\SG   f_n$ is Cauchy in $L^0$ and thus converges to some random variable $Y$. We define then $\SX  f:=Y$, thus extending the linear integral operator onto the whole $L^\phi$. By continuity, the Fourier transform of  $\SG f$ is given by \refrm{Ff}. 
\vspace{0pt}

 Conversely, if $Y=\lim_n \SG  f_n$ in probability then the convergence of the ch.fs.\ in view of \refrm{Ff} makes $f_n$ Cauchy in $L^\phi$, so by its completeness $f_n\to f$ in $L^\phi$. Again, \refrm{Ff} identifies the ch.f.\ of $Y$.
 \QED
  Such ``existential'' argument makes the stochastic integral quite enigmatic.
\subsection{Atomic approach}
The very L\'evy-Khinchin
representation of such  L\'evy process $\SX _t$   sheds some light on these formulas. We set aside a more general deterministic trend of bounded variation as well as a diffusion because they may overwhelm or even trash the effect of pure random jumps. Yet, we allow a possible delicate linear drift.
 Namely,
\[
\ln \E \ex{i\theta \SX _t}=t\left(i\delta \theta+\int_{\R} \left(e^{i\theta x}-1-i\,\theta\,\xx\right)\,\nu(dx)\right)
\]
where $\nu$ is an atomless L\'evy measure and $\xx$ is any measurable bounded function  such that  $\left|\xx-x\right|=O(x^2)$ near 0. Standardly, $\xx=x/(1+x^2)$, or $\xx=x\Ib{|x|\le 1}$, or $\xx=\signof x\,(|x|\wedge 1)$. Some special choices may fit specific distributions, though. For example, Zolotarev \cite[Intr.\ Thm.\ A]{Zol}) broke these standards  with a nonstandard $\xx=\sin x$ which entailed a commonly accepted new standard representation of the stable ch.f.  Variations that are caused by choices of $\xx$ can be simply incorporated in the shift parameter $\delta$. It happens for the Gamma process because $\xx$, being $d\nu$-integrable itself,  can be replaced just by 0.
\vv
We will consider shifts separately. The L\'evy measure of the Gamma process can be expressed in terms of the exponential-integral function
\be\label{E1}
E_1(v)=\int_v^\infty e^{-x}\,\frac{dx}{x}=\int_1 ^\infty e^{-vx} \,\frac{dx}{x}
\ee
with the density $e_1(x)=e^{-x}/x,\,x>0$.
\vv
The function has appeared independently in a variety of science applications. For example, in \cite{Pec} a suitably scaled  increment $E_1(a\,v^2)-E_1( a\,v^2_\infty)$ describes  the motion of a single meteoroid with reference to its atmospheric and preatmospheric velocities $v$ and $v_\infty$ (an algorithm for the inverse function was crucial).
\vv
Indeed, for $\Phi_1(|f|)<\infty$, the identities confirm formulas \refrm{Lf}  and \refrm{Ff}
\[
\begin{array}{rcl}
\Int_0^\infty \left(1-e^{-\theta u}\right)\,e_1(u)\,du
&=&\ln(1+\theta)\\
\Int_0^\infty \left(1-e^{\imath\theta u}\right)\,e_1(u)\,du
&=&
\frac{1}{2}\ln(1+\theta^2)-i\,\arctan( \theta)
\end{array}
\]
(the latter terms is just $\ln(1-\imath\,\theta)$).
Still, at this stage the existence of the integral process requires the consistency theorem.
\vv
In contrast, we can utilize a deterministic calculus of jumps (cf.\ App.\ \ref{jumps}). The integral (random, of course) emerges directly even on a quite arbitrary measure space  and the random measure $\SG \I A$ or process emerge unequivocally merely as by-products rather than an essence.
\vv
We still employ the arrivals  $(S_n)$ of a unit intensity Poisson process on $[0,\infty)$. Let $U_n$ be independent copies of random element $U$   with values in $(\BL,\ca L,\lambda)$ with a strictly positive density $p(u)$.  Also, let $(U_n)$  be independent of $(S_n)$. For example, partitioning $\BL$ into the union of disjoint  subsets $T_n$ with $\lambda T_n=1$, and introducing the probability measure
\be\label{pn}
\mu=\sum_n \lambda (T_n\cap \cdot) p_n, \quad\text{where}\quad  \sum_n p_n=1,
\ee
we define $U$ as the identity from $(\BL, \ca L,\mu)$ to $(\BL, \ca L,\lambda)$.
\vv
The inverse function  $H=E_1^{-1}$ entails the atoms $H_n$ and then the ``$H$-series'':
\be\label{H-series}
{\SG}f  =\sum_n H_n f(U_n),
\quad\text{where}\quad
H_n=H\left(S_n \,p(U_n)\right),
\ee
assumming the a.s.\ summability. The integrand-indexed process emerges first and only then it entails the stochastic process  and  the random measure as by-products:
\be\label {H-series paths}
\SG_t=\sum_n H_n \Ib{U_n\le t},\quad  \SG A=\sum_n H_n \Ib{U_n\in A}.
\ee
Straightforward computations  confirm \refrm{Lf} and \refrm{Ff} even for a general $(\BL,\ca L,\lambda)$.
\vspace{8pt}

Based on \cite{KalS}, we delineate a justification for the Gamma process on $[0,1]$, extended later to its quadratic variation. 
The display is more transparent than a similar one for the domain  $[0,\infty)$ which is just more elaborate. (One may need to switch to the Fourier transform if necessary.)
Indeed, let $U$ be a uniform random variable on $[0,1]$. By Fubini's theorem
\[
\E e^{-\SG f}=\E_S \prod_n \E _U e^{-H(\SG_n)f(U)}.
\]
Defining the function $\psi(x)$ by the formula
\[
e^{-\psi(x)}=\E e^{-H(x)f(U)},
\]
we can use the standard Poisson process $N_t$:
\[
\E e^{-\SG f}= \E e^{-\sum_n \psi(\SG_n)}=
\E e^{-N \psi}=\ex{-\int_0^\infty\left(1-e^{-\psi(x)}\right)\,dx}.
\]
That is,
\be\label{substH}
\E e^{-\SG f}= \ex {-\E\int_0^\infty\left(1-  e^{-H(x)f(U)}\right)dx}.
\ee
\begin{note}\label{arg}
The argument works for any positive pure jump L\'evy process. 
\QED
\end{note}
Specifically, for the Gamma process with $H$ being the inverse of $\nu(x,\infty)$ having the density $e^{-x}/x$,
the substitution $y=H(x)$ or $x=\nu(y,\infty)$ in \refrm{substH} yields
\[
\E e^{-\SG f}=\ex{-\lambda{ \ln(1+f)}}.
\]
Now, the compounding, captured earlier by the transform \refrm{FfK}, receives the a.s.\ representation for any admissible measure space $\BL,\ca L,\lambda)$:
\be\label{SSf}
\SG^{(\kappa)}f\df\sum_n K_n\,H_n f(U_n)
\ee
with independent  copies $K_n$ of a random variable $K$
(cf.\ App.\ \refrm{xk}).
\vv
With $H(x)=x^{-1/\alpha}$, the formulas are known as the LePage representations of an $\alpha$-stable integral or process \cite{LeP} (positive for $\alpha<1$ or compounded with random signs in the symmetric case). The H-series \refrm{H-series} for a general pure jump L\'evy process on an arbitrary continuous $\sigma$-finite measure space appeared in \cite{Ros} where even earlier references can be tracked down.
\vspace{5pt}

Alternatively but specifically on $\BL=\R_+$, we may first define $\SG_t$ on $[0,1]$ using \refrm{H-series} with $p=\I {[0,1]}$;  just using independent copies of a uniform random variable $U$:
next, with independent copies ${\SG}^{(n)}_t, t\in [0,1]$ at hand,  we may define inductively
\[
{\SG}_n \df {\SG}_{n-1}+{\SG}^{(n)}_1,\quad\text{then}\quad
{\SG}_{n+t} \df {\SG}_n+{\SG}^{(n+1)}_t, \quad n=1,2,...,\quad t\in (0,1].
\]
Now,  in view of Note \ref{arg} and a comment below \refrm{xk}, we look at the quadratic variation $[\SG] f=\sum_n H^2(S_n) f^2(U_n)=[\symS]f$. The proposed probabilistic form may offer more translucency than a coarse special function.
\begin{theorem}
The integrable function with respect to the  quadratic variation $[\SG]$ on $T\subset \R^+$ (or any $\sigma$-finite measure space) form the modular space
\[
\left\{f: \E \int_T \ln (1+2S_{1/2}f^2(x) )\,dx<\infty\right\}.
\]
\end{theorem}
\Proof 
W.l.o.g.\ we may and do assume that $T=[0,1]$. 
In view of \refrm{substH} we need to evaluate the following integrals:
\be\label{lnZ}
\int_0^\infty \left(1-e^{-H^2(x)\theta^2}  \right)\,dx=\int_0^\infty \left(1-e^{-x^2\theta^2}\right) \,\Frac{e^{-x}} x\,dx.
\ee
In the associated Laplace transform a standard Gaussian random variable $Z$ clarifies the appearance since  $Z^2\D S_{1/2}$.
With $\theta\df f(U)$, 
\[
\begin{array}{rcl}
\Int_0^\infty \Frac{1-e^{-x^2}}{x}\,e^{-x/\theta} \,dx &=& \E \Int_0^\infty \Frac{1-e^{\imath \sqrt{2}\,x Z}}{x}\,e^{-x/\theta} \,dx\\
&=& \E \Int_0^\infty \Frac{1-\cos (\sqrt{2}\,x Z)}{x}\,e^{-x/\theta} \,dx \\
&=& \E \Int_0^\infty \Frac{1-\cos x}{x}\,e^{-sx} \,dx \qquad\left(\text{with }\,s=(\sqrt 2\,\theta Z)^{-1}\right)\\
&=& \Frac 1 2 \,\E \ln(1+s^{-2})= \Frac 1 2 \, \E \ln (1+2 \theta^2 Z^2),
\end{array}
\]
concluding the computations.
\QED
\subsection{Thorin GGC class}
Thorin's ``generalized gamma convolutions'' \cite{Tho}, or  GGC,  are members of the smallest class $T(\R_+)$ of probability distributions spanned by Gamma distributions and closed under convolution and weak limit.  In \cite{Bar-NMS} the concept was extended to $T(\R^d)$. In our context, the integrals ${\SG}f$ of positive simple functions serve as a generator of the Thorin class. Apparently, the class contains constants, obtained, e.g., from the Law of Large Numbers, which are not representable by integrals.
\vv
Setting degenerate distributions aside, the logarithm of the Laplace transform of a Thorin's GGC has the representation
\[
\ln \widetilde{\mu}(\theta)=-\int_0^\infty \left(1-e^{-\theta x}\right)\,\nu(dx),\quad \int_0^\infty (1\wedge x)\,\nu(dx)<\infty.
\]
Also, $\nu(dx)=\Frac{k(x)}{x}\,dx$ with a completely monotonic $k(x)$, i.e., for all nonnegative integers  $n$, $(-1)^n k^{(n)}(x)\ge 0$.  In particular, by Bernstein Theorem, $k(x)$ is the Laplace transform of a $\sigma$-finite measure on $\R_+$, called Thorin measure. The following equivalence is well known (e.g., \cite[Prop.\ 1.1]{JamRY} with references to earlier sources.)
\vv
\begin{theorem}
The probability distributions of Gamma entities ${\SG}_f$ coincide with nondegenerate Thorin GGC.
\end{theorem}
\Proof
By construction, the law of ${\SG}_f$ belongs to the Thorin class. Let us find its characteristics. W.l.o.g.\ we assume that $T=[0,1]$.
\vv
Replace $f$ by $\theta f$ in formula \refrm{substH}, then substitute $y=H(x)f(U)$, yielding
\[
\ln \E e^{-\theta\xi Lf}=-\int_0^\infty \left(1-e^{-\theta y}\right) \Frac{k(y)} {y}\,dy ,\quad\mbox{where} \quad
k(y)=\E \ex{-y/f(U)}.
\]
Clearly, $k(y)$ is completely monotonic and $k(y)=\widetilde{\sigma}$, where
$\sigma=\sigma_f$ the probability distribution of $1/f(U)$.
\vv
Conversely,  a Borel probability (w.l.o.g.) measure $\sigma$ with c.d.f.\ $G(x)$ defines the function  $f(u)=G^{-1}(u),\,u\in [0,1]$. Then, reversing the above steps, we infer that the law of ${\SG}_f$ is Thorin's GGC with Thorin measure $\sigma$.
\QED
The difference between an ``entity'' marked by a function or scalar and a ``stochastic process'' can be illustrated by the following examples.
\vv
{\bf Example 1}.
The radical of a positive Gamma process may serve as a random time replacement, although flawed because it that lacks independent stationary increments), in a scaled Wiener process; producing one-dimensional distributions:
\[
\ca L\left(\symS_t\right)=
\ca L\left(\sqrt{2}\,W\left(\SG_{t/2}^{1/2} \right)\right)
\]
for \underline{each single nonnegative $ t$}
but the processes are not equidistributed. In contrast, the independent (proper) subordinator $\SG_t$ yields the L\'evy process process $W\circ \SG$.
\[
\E e^{i\theta W(\SG_t)} =
\E\left((1+\theta^2Z^2/2)^{-t}\right),\quad  \text{where $Z\sim$   N(0,1)}.
\]
Its paths are almost everywhere discontinuous.
\QED
{\bf Example 2}. Let $\alpha\in (0,2)$ and $f(x)=x^{-1/\alpha}$, $x>0$.  Then the ``entity'' $\SG_f\df \SG f$ is a positive $\alpha$-stable random variable for $\alpha>1$ and $\symS_f\df \symS f$ is symmetric $\alpha$-stable for $\alpha<2$. Indeed, just compute the Laplace or Fourier transforms.
Yet, these ``frozen'' values of the Gamma integral processes have nothing in common with the L\'evy stable processes or integrals, except for these solitary coincidences.
\QED
A survey of the ample GGC class can be checked in \cite{JamRY} along with references therein. The reward  Gamma processes \refrm{SSf} appear also there under the name ``{\em subordinators
$\SG_t(G)$}'' (processes, not single entities);
 the only difference being $K = \frac 1 G$ (see \cite[(50) and next]{JamRY}).  Since only positive Gamma process and positive integrands are used, and we consider $G\ge 0$, we convert \refrm{FfK} to the Laplace transform and perform some calculus on $\R_+$, using just the logarithm and $t=1$;  with the probability distribution  $\nu$ of $G$:
\[
\E \int_0^\infty \ln\left(1+\frac {f(x)} G  \right)\,dx=
\int_0^\infty \int_0^\infty \ln\left(1+\frac {f(x)} u  \right)\,dx\,\nu(du)
\]
\section{Paths and integrability}
\subsection{Jumps}
While every L\'evy process has a rcll version,  the pointwise representation \refrm{H-series paths} is rcll by a deterministic argument (cf.\ Appendix \ref{jumps}).
However, no path-continuous version exists since ${\SG}_t-t$ is a martingale of finite variation with all finite $p$-moments, $p>0$.
 \begin{theorem}\label{Lip}
The paths of  $\SG_t$ are right-Lipschitz. That is, for every $a>0$ there is a random quantity $K_{a}$, although nonintegrable, such that
\[
0\le {\SG}_t-{\SG}_a\le K_{a} \,(t-a),\quad t>a.
\]
\end{theorem}
The following identity holds even in a more general context
 (cf., e.g. \cite[Lemma 3.3.1]{RosBI}).
 \begin{lemma}\label{abcf} 
 $\E\left[{\SG}_a f\,|\,{\SG}_c  \right]=
 \left(\lambda_a f\right)\,\Frac{{\SG}_c}{c}\quad 0\le a\le c.$
\end{lemma}
\Proof
Indeed, points $\ca P_n= \set{kc/n: k=1,\dots, n-1}$ generate a uniform partition of the interval $[0,c]$. The random variables ${\SG}_{kc/n}-{\SG}_{(k-1)c/n}$ are exchangeable conditionally, given ${\SG}_c$. Hence, for $a\in \ca P_n$,
\[
\E[{\SG}_a|{\SG}_c]=\frac{a}{c} \,{\SG}_c,
\]
extending to any $a$ by the right continuity, and then to integrals by routine.
\QED
{\bf \Proof of Theorem \ref{Lip}}. Consequently
for a fixed $a\ge 0$ we obtain the reverse martingale
\[
\eta_t=\frac{{\SG}_t-{\SG}_a}{t-a}, \quad t>a,
\]
 with respect to the descending filtration $\ca R_t=\sigma\set{{\SG}_b-{\SG}_a:b\ge t}$. It
converges a.s.\ as $t\to a$ since it is bounded in $L^1$.  By the Kolmogorov 0-1 Law, the limit is a non-random constant which is 0 due to its Laplace transform. Thus, $K_a\df\sup_{t>a} |\eta_t|<\infty$ a.s.\ but $\E K_a=\infty$  since otherwise the reverse martingale $\eta_t$ would be closed but it is not.
\QED
Let us consider now the inverse process $R\df {\SG}^{-1}$, cf.\ \refrm{app:inv}, with paths illustrated by truncated series as in Section \ref{cut-off}. The cut-off already suggests the continuity of paths and a sort of possible ``explosion'' which prompts the computation of suitable measures such as the hyper-exponential moments
$\E R_t^{\theta R_t}$, indeed exhibiting a sort of growth phase transition at $t=1/e$ from ``fast'' to ``hyper-fast''.

\begin{figure}
\begin{center}
  \includegraphics
  [width=0.8\textwidth]
  {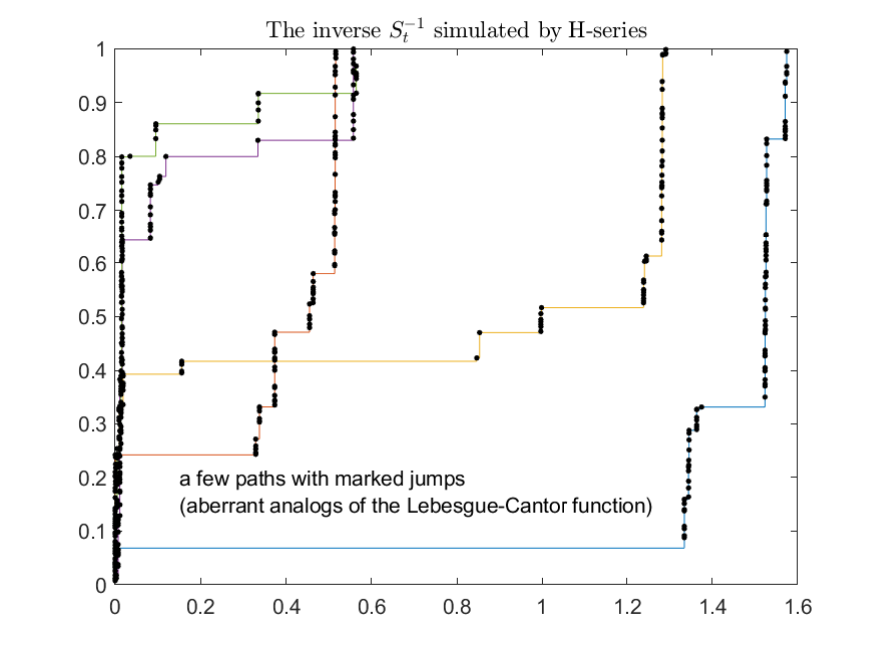}
  \caption{The range is truncated by 1.}
  \end{center}
\end{figure}
\vv

\begin{theorem} Let $R={\SG}^{-1}$.
\begin{enumerate}
\item Almost all paths of $R_t$ are continuous.
\item
$\E R_t^{\theta R_t}<\infty$ for $\theta<1$ and $t>0$; 
so, the exponential moments $\E e^{\theta R_t}$ exist.
\item
$\E R_t^{R_t}<\infty$ iff $t<\frac 1 e$.

\end{enumerate}

\end{theorem}
\Proof
For the continuity see Appendix Note \ref{app:inv:cont}. \vv

Next, applying the identity $\E \phi(Y)=\phi(0)+\int_0^\infty \phi'(x)\,\PR(Y>x)\,dx$ to the function $\phi(x)=x^{\theta x}$ with $0^0\df 1$,
\[
\E R^{\theta R}-1=\theta \int_0^\infty (1+\ln x)\, x^{\theta x}\,\PR(R>x)\,dx=
\theta \int_0^\infty (1+\ln x) x^x\PR({\SG}(x)\le t)\,dx
\]
It suffices to truncate the integral using an arbitrary large lower limit, i.e., integrating over $[c,\infty)$.
For a fixed $f$ we also have
\[
\PR({\SG}(x)\le t)=\frac{1}{\Gamma(x)}\int_0^t u^{x-1} e^{-u}\,du\approx \frac{t^x}{\Gamma(1+x)}.
\]
Thus, by Stirling's formula,
 \[
\E R^{\theta R} \approx \int_c^\infty  \ex{(\theta-1)\, x\ln x+x\,(\ln t+1) -\frac{1}{2}\,\ln x+\ln\ln x}   \,dx,
 \]
proving the last two statements.
\QED

\subsection{Moments}
As usual, $\|X\|_p$ or $\|f\|_p$ denote the $p$-norms or $p$-F-norms (when $p<1$).
\vv
From \refrm{Ff} we infer that $X=\SG f \Ib{|f|\le c}$ possesses the exponential moments $\E e^{\theta X}$ whenever $|\theta|<1/c$ for any $c>0$, so all moments $\E |X|^p$  exist. On the other hand, $\lambda\{|f|>c\}<\infty$. Therefore, in studying the moments of $\SG f$ it suffices to confine to $f$ with a support $T$ of finite measure. Furthermore,  we may have to relinquish some hard analysis (that is, close estimates of suitable constants). The connection between the differentiability (quite straightforward in our case) of the characteristic function and moments is well known (cf. \cite[Sect.\ 2.3]{Luk}) as well as are the difficulties related to odd absolute moments; and all the more, to general non-integer $p$-moments.
\vv
For illustration, we observe that   $\E({\SG}f)=\lambda f$ for $f\ge 0$ and it stands to reason to believe that the identity holds for a general $f$.  However, the mere existence of the mean $\E|\SG f|$ requires a  proof of integrability of $f$. Although the ch.f.\ of $|\SG g|$ can be expressed \cite[(2.3.11)]{Luk} in terms of the derivative $F'$ of the ch.f.\ of $\SG f$, but a derivation  of the integrability, $\lambda|f|<\infty$, seems very difficult. We will show details in a more general context.

\begin{theorem} \label{pint}
Let $p>0$ and $\SG f$ exist for $f\in L^0(\BL,\ca L,\lambda)$.

Let $|\supp f|<\infty$. Then there exist positive constants $c, C$ such that
\[
c||f||_p \le ||{\SG} f||_p \le C ||f||_p.
\]
Consequently,   $\E |{\SG}f|^p<\infty$ if and only if for  any or some $c,c'>0$,
\be\label{Sfp}
\mbox{(a) }\lambda |f| \,\Ib{|f|\le c}<\infty\quad\mbox{and}\quad
\mbox{(b) }\lambda |f|^p\Ib{|f|> c'}<\infty.
\ee
\end{theorem}

We need a few auxiliary results, all under the assumptions of the theorem.
\begin{lemma} \label{pint1}
 $\E\left|\SG f\Ib{|f|\le c}\right|^p<\infty$ for every $p>0$.  In addition,
\[
\E\left(|\SG f|^p\right)<\infty \quad \iff \quad \E\left((\SG f^2)^{p/2}\right)<\infty.
\]
Further, by Baire's category argument, there exist positive constants $c_p, C_p$ such that for every $f$,
\be\label{Bcat}
c_p \left\|f(\SG f^2)^{1/2}\right\|_p\le \|\SG f\|_p\le C_p \left\|(\SG f^2)^{1/2}\right\|_p.
\ee
\end{lemma}
\Proof The first claim stems from the exponential moment
$\E \ex{\theta \SG f\I{|f|\le c}}$ that exists when $|\theta|<1/c$.
\vv
In general,  $\|X-X'\|_p\le 2\|X\|_p$ for a random variable $X$  and its independent copy $X'$.  On the other hand, if $\E|X-X'|^p<\infty$ then by Fubini's theorem $\E|X-x'|^p<\infty $ for some number $x'$. That is,
$\|X\|_p\le \|X-x'\|_p+|x'|$.
\vv
Consequently, the moments $\E|\SG f|^p$ and $\E|\symS f|^p$ exist simultaneously (in contrast to the mere existence of the integrals).
Therefore, we may employ the compounded series \refrm{FfK} with the symmetric i.i.d.\ Rademacher random variables (i.e., with the symmetric $(\pm 1)$-Bernoulli distribution and in virtue of Fubini's theorem focus just on the Rademacher series. All their moments are comparable, in particular, to the second moment. Translating this back to the language of Gamma integrals, the second statement has been proved.
\vv
Finally, the category argument is valid since  we deal with at least F-spaces $(p<1$) if not Banach spaces $(p\ge 1$).
\QED
The following  well known formula, e.g., \cite[(2.4.1) and (2.4.4)]{Luk}, needs just a little adaptation to our context. Still, a direct derivation shows how the symmetric and skew parts interfere and complete each other.
\begin{lemma} Let $p\in\N$.
If the integrals $m_\ell=\lambda f^\ell$ exist for $\ell=1,...,n$ then
 the $p$\tss{th} moment
\be\label{mun}
\alpha_p\df \E\left(|\SG f|^p\right) =p!\Sum_{k=1}^p \Sum_{\wek j}
\Frac 1{\wek j\,!}
 \displaystyle \prod_{\ell=1}^{p-k+1} \left( \Frac {m_\ell}{\ell} \right)^{j_\ell},
\ee
with the integer vector indexes $\wek j=(j_1,\dots, j_{p-k+1})$, $\wek j\,! \df j_1!\cdots j_{p-k+1}!$ and subject to the constraints
\be\label{j}
\Sum_{\ell=1}^{n-k+1} j_\ell=k \quad\text{and}\quad
\Sum_{\ell=1}^{n-k+1}\ell  j_\ell=p.
\ee
\end{lemma}
\Proof As usual,
$ \alpha_p :=\E (\SG f)^p = \imath^{-p}\, F^{(p)}$, where
$F(\theta)=\E e^{\imath \theta X}=e^{\psi(\theta)}$ and
$
\psi(\theta)=
-\frac{1}{2}\tau \ln(1+ f^2\theta^2)+\imath\,\tau  \arctan(f\theta)
$.
\vv
The first derivative involves a single complex valued function
\be\label{1}
\psi'(\theta)=\lambda f(c(f\theta))\quad\text{with}\quad
F'(0)=\psi'(0)=\imath\, \lambda f
\ee
because
\[
\psi'(\theta)=
 -\lambda \frac {f^2 \theta}{1+ f^2\theta^2}+\imath
 \lambda \frac {f}{1+ f^2\theta^2}=
\lambda f \Big(-a(f\theta) +\imath b(f\theta)\Big),
\]
yielding  $c=-a+\imath b$ with
\[
b(\theta)=\Frac 1{1+\theta^2}\quad\text{and}\quad
a(\theta)=\theta\,a(\theta).
\]
Therefore, for $\ell=1,...,p$, we have $\imath^{\ell+1}\, c^{(\ell)}(0)=(-1)^{\ell+1}\ell!$  because
\[
b^{(\ell)}(0)
=\defifsmall {(-1)^m \ell! } {p=2m} 0 {\ell=2m+1}\,\,
a^{(p)}(0)
=\defifsmall {0} {\ell=2m}  {(-1)^\ell ! } {\ell=2\ell+1}
\]
Hence,  and  from \refrm{1} we infer by the Chain Rule that
\be\label{chain}
 x_\ell\df \psi^{(\ell)}(0)=(-1)^\ell\frac{(\ell-1)!}{\imath^{\ell}}\,m_\ell,
 \quad \ell = 1,..., p.
\ee
Using Bell's polynomials $B_{p,k}$, $1\le k\le p$,
\[
B_{p,k}(x_1,\dots,x_{p-k+1})=n!\sum_{\wek j} \Frac 1{\wek j\,!} \prod_{\ell=1}^{p-k+1} \left( \Frac {x_\ell}{\ell!} \right)^{j_\ell}
\]
that appear in Fa\`a di Bruno's formula
\[
F^{(p)}=F
\Sum_{k=1}^p B_{p,k}\Big(\psi',\psi'',...,\psi^{(p-k+1)}\Big)
\]
with quantities \refrm{chain}, we obtain\refrm{mun}.
\QED
\begin{corl}\label{pint2}
If $p\in \N$ and $f\ge 0$, denoting $\mu_p^*=\Max_{1\le \ell \le p} \|f\|_\ell$,
\be\label{pnorm}
\mu_p^*\le \|\SG f\|_p \le (p!)^{1/p} \,\mu_p^*.
\ee
\end{corl}
Indeed,  the right constant follows from the $p$-homogeneity of the Bernoulli polynomials. After $\mu_p^*$ is factored out, the remaining quantity is just the $p$th moment of an exponential $S_1$, or $p!$.
\vv
On the other hand, for each integer $q\le p$, dropping all terms except one for $k=1$, the restriction holds only for $j_q=1$ and $j_l=0$ for $\ell\neq q$. So,
\[
\alpha_q \ge q! \frac{m_q}{q}=(q-1)! m_k.
\]
That is, $\|S f\|_p \ge \Max_{ q\le p} \left(((q-1)!)^{1/q} \|f\|_q\right)\ge \mu_p^*$.
\QED
{\bf Proof of Theorem \ref{pint}.}~ By Lemma \ref{pint1} we may confine to $f\ge 0$ with $\lambda \supp f=c<\infty$, so $\mu_p^*\sim \|f\|_p$. Hence, for an integer $p$, Corollary \ref{pint2} provides the sought-for isomorphism of both $L^p$ spaces.
For a non-integer $p$ we apply the Riesz-Thorin interpolation theorem, so the isomorphism holds for every $p\ge 1$.
\vv
Let $p<1$. We  choose a positive $p$-stable random variable $S$, normalized to secure the identity $\E e^{-\theta S }=e^{-\theta^p}>$, $\theta\ge 0$, and independent of the Gamma process ${\SG}_t$. Then, by Fubini's Theorem,
\[
L(\theta)=\E e^{-\theta ({\SG} f)^p} =\E e^{-s\theta^{1/p} S \,{\SG}f}=\E \ex{-\int_T \ln\left(1+s^{1/p} S \,f\right)}=
\E M.
\]
where the random variable $M=M(S ,\theta)\le 1$. We compute and estimate the derivative:
\[
-L'(\theta)=\frac{1}{p\theta}\,\E \left[M(S,\theta)\, \lambda \frac{\theta^{1/p} S  f}{1+\theta^{1/p} S f}  \right]\le
\frac{1}{p\theta}\,\E \left[\lambda \frac{\theta^{1/p} S f}{1+\theta^{1/p} S  f}  \right].
\]
Since
$\frac{u}{1+u}\le 1-e^{-u},\,u\ge 0$,  therefore by Fubini's Theorem and then, using the inequality
$1-e^{-u}\le u,\, u\ge 0$, we continue:
\[
-L'(\theta)\le \frac{1}{p\theta}\,\lambda\E \left[1-e^{-\theta^{1/p} S  f}  \right]=
\frac{1}{p\theta}\,\lambda \left[1-e^{-\theta f^p}  \right]\le
\frac{1}{p}\,\lambda  f^p
\]
Letting $\theta\to 0$, we obtain the estimate
\[
\E ({\SG} f)^p\le \frac{\lambda f^p}{p}.
\]
Since $\lambda \supp f=c<\infty$, then by stationarity we may assume that $\supp f=[0,c]$. So, by Fubini's Theorem and Jensen's inequality on a (randomized) probability space
\[
\E({\SG}f)^p=\E \left(\int_0^c f\,d{\SG}\right)^p= \E \,{\SG}^p_c \left(\int_0^c f\,\frac{d{\SG}}{{\SG}_c}\right)^p\ge
\E \,{\SG}^{p-1}_c \int_0^c f^p\,{d{\SG}}
\]
The latter expectation is equal to
\[
\E\left( \E\left[ \,{\SG}^{p-1}_c \lambda_c f^p \,d{\SG}\Big|{\SG}_c\right]\,\right)=
\frac 1 c\,  \E {\SG}^{p}_c \lambda_c f^p=\Frac{\Gamma(p+c)}{c\Gamma(c)} \lambda{f^p}
\]
in virtue of Lemma \ref{abcf}.
\QED
\begin{remth} [Discussion of constants]{~}

\begin{enumerate}
\item {\em Let $\lambda \supp f= 1$. Then $\mu_p^*=\|f\|_p$.  If $f\ge 0$, then the  left constant in \refrm{pnorm} can be improved to $((p-1)!)^{1/p}$.}

\begin{quote}
{\small Indeed, the inequality in the last line of the Corollary's proof now involves the maximum of the product of two increasing sequences, hence the maximum is attained at the product of the maxima.
\vv
Further, these are the best constants. The right one is attained for $f=\I T$ wit $\lambda T=1$. In regard to the lower bound, there is a sequence of functions $f_j$ such that $\lambda f_k^n\to 1$ but $\lambda f_k^\ell\to 0$ for $\ell<n$. For example, simplifying the context to $\BL=\R_+$.
\[
f_k(x) = \Frac 1 {x^{1/n}\ln ^{1/n}(k)} , \quad \frac 1 k\le x\le 1.
\]}
\end{quote}
\item {\em Still with $\lambda \supp f=1$, an even power $p$, and $f$ admitting negative values but with $\lambda f=0$, the constants obtained above stay.}

    \begin{quote}
   {\small  The upper bound relies just on the triangle inequality. For the lower bound, with $m_1=0$ and even $p$, formula \refrm{mun} contains only even powers of odd moments, due to the second constraint. Thus, all summands are nonnegative, validating the ``dropping argument'' in the proof of Corollary \ref{pint2}.}
    \end{quote}

\item {\em Odd $p$ and $f$ with possible negative values face a classical obstacle. We found no decent hard estimates.}
\item
{\em Let  $f\ge 0$ with a support of positive finite measure $\lambda T=c$ or, equivalently, let $\SG$ have intensity $c$. Then $\|\SG f\|_p\approx \frac p e \|f\|_p$, or more precisely,
\be\label{Gfn:bounds}
\frac{\Gamma(c+p-1)}{\Gamma(c)}\, \lambda f^p \le \E \left(\SG f^p\right)\le \frac{\Gamma(c+p)}{\Gamma(c)}\,\lambda f^p,
\ee
with the best constants.}

\begin{quote}
{\small By rescaling $\lambda\mapsto c\lambda$ results in the replacement $m_\ell \mapsto c m_\ell$ in the moment formula \refrm{mun}, which is tantamount to the change of intensity $1 \mapsto c$, or to the new $T=c T_c$ where $\lambda T_c=1$. Again, if $m_p\le 1$ then $m_\ell\le m_p\le 1$, so
\[
\begin{array}{rcl}
\E\left( (\SG f)^p\right)&=&p!\Sum_{k=1}^p c^k \Sum_{\wek j}
\Frac 1{\wek j\,!}
 \displaystyle \prod_{\ell=1}^{p-k+1}
  \left(\Frac {m_\ell}{\ell}\right) ^{j_\ell}\\
  &=& n!\Sum_{k=1}^p c^k \Sum_{\wek j}
\Frac 1{\wek j\,!}
 \displaystyle \prod_{\ell=1}^{p-k+1}
  \Frac {1}{\ell^{j_\ell}}=\E\left((\SG \I {T_c})^p\right)=\frac{\Gamma(p+c)}{\Gamma(c)}.
  \end{array}
 \]
 In particular,  the $p$\tss{th} moment of the $\Gamma(c,1)$ distribution yields \refrm{Gfn:bounds}. Stirling's formula ensures the equivalence of the $p$\tss{th} norms.}
 \end{quote}

\end{enumerate}
\end{remth}
\subsection{Simulation}\label{cut-off}
\subsubsection{By partitions}
A common approach goes back to Wiener's concept \cite{Wie} of a stochastic integral  $\SX f:=\int f\, d\SX$ on $\R_+$ that mimics approximation by step functions, constant between partition  points $\{0=u_0<u_1<\cdots <u_{n-1}<u_n=1\}$ of, say, $[0,1]$:
\[
 g=\sum_{j=1}^n a_j \I{\left(u_{j-1},u_j \right ]}
 \quad\text{yielding}\quad  \SX g:=\sum_{j=1}^n a_j
 \Big( \SX_{u_j}-\SX_{u_{j-1}}  \Big).
\]
Typically, the values at the left end points are used, e.g., $a_j=f(u_{j-1})$, but not necessarily (see \refrm{aj} below).
For a second order process $\SX$, the variance $\Var(\SX g-\SX g')=\|g-g'\|_2^2 $ yields $\SX f$ by completion. Consequently, the goodness of approximation can be measured by the variance
\[
\begin{array}{rcl}
\varepsilon_n =\Var(\SX f-\SX g)= \|f-g\|_2^2=\Sum_{j=1}^n \Int_{u_{j-1}}^ {u_{j}}
\left|f(t)-f\left(u_{j-1}\right)\right|^2\,dt.
\end{array}
\]
The best least squares approximation of $f\in L^2[0,1]$ by a function with finitely many values that span a finite field $\ca T$ is given by the conditional expectation, a.k.a.\ projection onto $L^2(\ca T)$,
\[
\E\big|f-\E[ f|\ca T]\big|^2=\inf_{g\in L^2(\ca T)} \E|f-g|^2.
\]
In this case the step function $g=\E[f|\ca T]$ possesses the values
\be\label{aj}
a_j =\frac 1 {u_j-u_{j-1}}\int_{u_{j-1}}^{u_j} f(x)\,dx.
\ee
For a L\'evy process, in virtue of its stationarity one may consider a monotonic rearrangement of $f$, i.e., w.l.o.g.\ one may assume that $f$ is monotonic. Then $a_j$ is contained between the values at the end points.  A uniform partition $u_j=j/n$ is most convenient for a L\'evy process since the increments are i.i.d.\ random variables.
\vv
A pure jump L\'evy process exhibits its own regime of jumps, contradicting the simulation by the above method which  places incorrect jumps at wrong positions in an artificial manner: true jumps are not Gamma variables and occur elsewhere.
 Still, such discrepancy is beyond detection of a naked human eye when the resolution is high enough. We only have a ``soft'' qualitative confirmation: for a nonnegative $f\in L^2$ and step functions, described above, $\SG f_k\to \SG f$ a.s.\ and in $L^2$.
At the same time, the question about the quantitative error $\varepsilon_n$ seems to have rather elusive answers. 
\vv
Indeed, when the given $f$ itself resembles a step function, the error can be arbitrarily small, even 0 eventually. Yet, for a typical ``orderly'' function such as a power $x^p$ with  $p>0$ or the exponential $e^{cx}$ the nullity is of order $O(n^{-2})$, rather weak. To wit, if $f\in C^1[0,1]$ with $f'\in L^2[0,1]$, using suitable intermediate points $\xi_j\in \left(\frac {j-1} n, \frac j n\right] $,

\hspace{85pt}
$
\varepsilon_n  = \Frac 1 {n^3}\, \Sum_{j=1}^n
\left|f'\left(\xi_j\right)\right|^2
\approx \Frac{\|f'\|_2^2}{n^2}.
$
\vv
Further, on the edge of square integrability the order of nullity may worsen significantly. For example, for $f(x)=x^{-p}$, where $2p>1$, the error $\varepsilon_n \sim n^{2p-1}$.
\subsubsection{By H-series}
In contrast,  the truncation ${_{_N}}\SG_t$ of  series \refrm{H-series} to the finite sum for $n\le N$ yields the easily controllable remainder
$R_N=(\SG - {_{_N}}\SG_t)f$, subject to moment assumptions. The assumption below reduces to $f\in L^p$ for any $p>0$ if $f$ has a support of finite measure.
\begin{rem}
The more slowly series \refrm{pn} converges, the better (faster) the approximation. In other words, the more closely the distribution of chosen $U$ on the half line resembles the uniform distribution, the better simulation. So, Pareto $U$ with the density $(p+1) x^{-p}\I{x\ge 1}$, where $p>1$  just barely, is more efficient than an exponential $U$ with the density $e^{-x}$.
\end{rem}
\begin{theorem} \label{rnp}
As $N\to\infty$, $R_N\to 0$ a.s.  Further, for an even integer $p$,  if a nonnegative $f\in L^q$ for $q=1,...,p$  and $M_p=\Max_q \lambda f^q$ then
\[
\E R_N^p \le \frac{C_p M_p}{(p+1)^N},
\]
for some constant $C_p>0$, i.e., the $L^p$-error tends to 0 geometrically.
\end{theorem}
\Proof
The a.s.\ convergence is obvious. Then, for Poisson arrivals  $(S'_n)$ that are independent of $(S_n)$, we have
\[
R_N\D \sum_n H(S'_N+S_n) f(U_n).
\]
Next, we note the inequality, clearly true for $x=0$:
\be\label {x+y}
H(x+y)\le e^{-x} H(y).
\ee
 Indeed, let $x>0$ and  put $u=H(y)$, i.e., $y=E_1(u)$, and also $c=e^{-x}<1$. The inequality is equivalent to
\[
-\ln c \ge E_1(cu)-E_1(u)
\]
The $u$-function on the right is decreasing hence its limit at 0 becomes the supremum. Based on the representation (cf.\ \cite[5.1.11]{AbrS})
\[
E_1(u)=-\gamma -\ln (u) -\sum_{k=1}^\infty \frac{(-1)^k} {k\,k!}\,u^k
\]
where $\gamma$ is the Euler constant, that limit equals $-\ln c$, thus proving \refrm{x+y}.
\vv
Hence, by Fubini's Theorem,
\[
\E R_N^p \le \E e^{-pS_N}\,\E \SG ^p f,
\]
establishing the estimate in virtue of the upper bound in Theorem \ref{pint}.
\QED

\begin{remth}{~}\sf

\begin{enumerate}
\item {\em  In general, inequality \refrm{x+y} may have no useful analog}.
\begin{quote}
{\small  E.g., it fails for an $\alpha$-stable process.  Nevertheless, there may  exist $N$ (and it does in the case of a stable process)  such that $H(S_n)$ is integrable for $n>N$, so the remainder $\E R_N\to 0$.}
\end{quote}
\item
{\em  Still, possibly, none of $H(S_n)$ may belong to any $L^p,\, p>0$.}
\begin{quote}
{\small For example, the L\'evy density $\approx \ex{-\ln^{1/q}(x)}$, $q>1$, at $\infty$ yields $H(x)\approx \ex{-\ln^q(x)}$ near 0.}
\end{quote}
\end{enumerate}

\includegraphics
[trim=1in 3in 1in 3in,clip, width=\textwidth]
{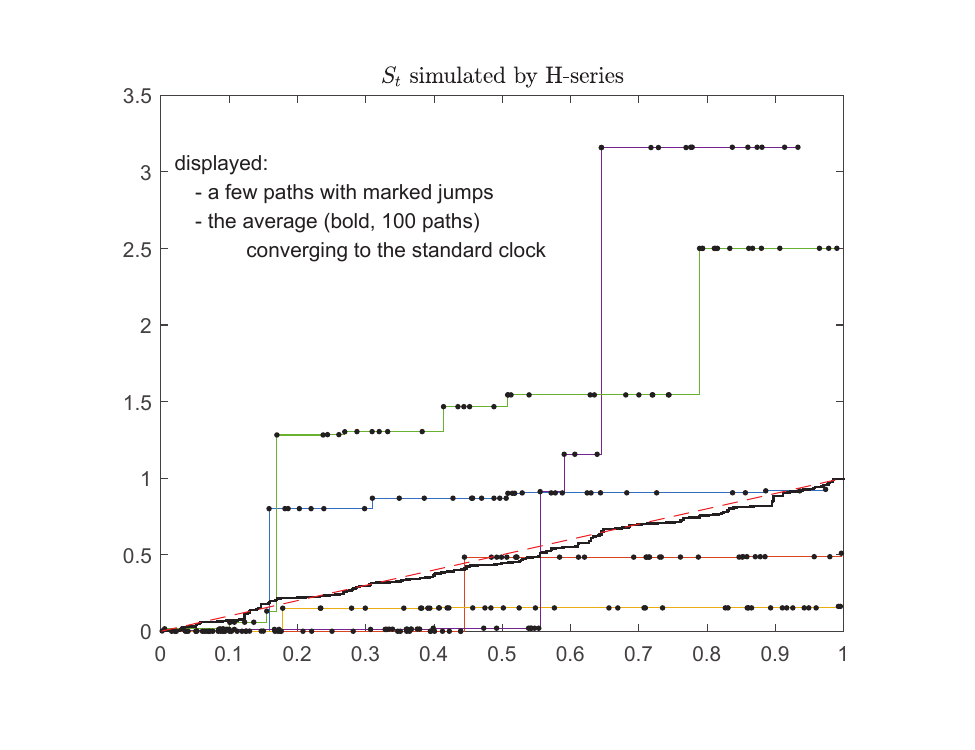}

In the above simulation we use a sufficient cut-off of the series.
In old days, the inversion of the exponential integral function was quite  challenging (cf., e.g., \cite{Pec}). While the modern software such as Matlab provides the ``{\tt expint}'' function, we found no built-in efficient inverse. However, a down-to-earth  ``halve-and-check'' algorithm works surprisingly fast, being as accurate as desired. Note that in the cut-off jumps points are order statistics and jump heights are not Gamma variables.

\end{remth}
\subsection{Linear trend}

An added deterministic trend $Y_t=\SX _t+b_t$ of bounded variation that ensures the infinite divisibility has an impact on the integrands, even for a general L\'evy process. The mere existence in $L^0$ of the integral $Yf$ implies the existence of $\delta f=\int_T f\,db$, cf.\ \cite[Th.\ 8.3.1]{KwaW} with a quite formidable proof. In contrast, in our case the argument is rather simple. Indeed, for $f\ge 0$ the imaginary part yields the modular space $L^{\phi_0}\cap L^{db}$, in particular securing the existence of the integral $\lambda f$.
\vv
\begin{rem}\label{f2}{~}

\begin{enumerate}
\item Consider a trend-free $\SX $.
\begin{enumerate}
   \item If  $\beta=0$ then
   the analog of Theorem \ref{pint} for $\symS f$ holds for  $f\in L^{\phi_2}$ (not necessarily positive) with condition (a) in \refrm{Sfp} replaced by the condition ``\,($\wt{\text{a}}$) $\lambda f^2\Ib{|f|\le c}<\infty$''.
   \begin{quote}
   {\small \sf Indeed, for every $p>0$ by Khinchin's inequality,
   \be\label{symp}
   c_p||({\SG}f^2)^{1/2}||_p\le ||\symS f||_p\le C_p ||({\SG}f^2)^{1/2}||_p,
   \ee
   for some constants $c_p,C_p$, independent of $f$.}
   \QED
   \end{quote}
   \item Let $\beta\neq 0$ and $f\in L^{\phi_1}$. Then $\SX f$ and $f$ are simultaneously $p$ integrable for every $p>0$.
       \begin{quote}
      {\small \sf Indeed, let $\SX  f$ exist. Then $f\in L^{\phi_0}$, which implies condition (a) of \refrm{Sfp} of Theorem \ref{pint} (stronger than ($\wt{\text{a}}$) above).  The equivalence follows from the previous statement.
      \QED
       }
       \end{quote}

   \end{enumerate}
\item Consider a nontrivial trend $b_t$.
\begin{enumerate}
\item
For a linear trend $b_t=\delta t$, the $p$-integrability occurs simultaneously for $Yf$, $\SX f$, $\SG f$,  and $f$.
\begin{quote}
{\small\sf Indeed, repeat the above argument.}
\QED
\end{quote}
\item {\sf A general trend $b_t$ with bounded variation entails the space $L^1(db)$ which interferes in various ways with the underlying modular spaces $L^{\phi_1}, L^{\phi_2}$, and $L^{\phi_0}$ through condition (a). We omit the discussion.
    }\QED
\end{enumerate}
\end{enumerate}
\end{rem}
\section{Gamma martingale fields}
A L\'evy process $\SX _t$ entails the exponential (pathwise for a PJP, cf.\ Appendix \ref{jumps}):
\be\label{expx}
{\ca E}_{\xth} \SX _t\df a_t(\xth) e^{-\xth \SX _t},
\ee
where $a_t(\xth)=e^{t\psi(\xth)}$ and $\psi(\xth)=-\ln \E e^{-\xth \SX _1}$.
Besides the self-explanatory algebraic form, independent $\SX_t$ and $\SY_t$ satisfy the identity
\[
{\ca E}_{\xth} \SX _t\,\cdot\,{\ca E}_{\xth} \SY _t\D
{\ca E}_{\xth} (\SX _t + \SY_t)
\]
(with possible versions that satisfy the equality pointwise on a product probability space).
Subject to the integrability assumption, hopefully beyond the trivial $\xth=0$,  the exponential becomes tautologically  a unit mean positive  martingale (or rather a martingale field)  w.r.t.\ the natural filtration $\ca G_t$.
For the Gamma process:
\[
\begin{array}{rcll}
\vspandexsmall
M_t(\xth)&=& {\ca E}_{\xth} \SG _t
& \text{where }\, a_t(\xth)=(1+\xth)^t
\,\text{ and }\, \xth>-1;\\
\wt{M}_t(\xth)&=&{\ca E}_{\xth} \symS _t
& \text{where }\,\wt{a}_t= (1-\xth^2)^{t/2}
\,\text{ and }\, |\xth|<1.
\end{array}
\]
With a linear trend, ${\ca E}_{\xth}(\SX _t-\delta t)$ is a martingale with $a(t)=\ex{t\psi(\theta)+\delta t}$ iff $\psi(\xth)= \delta \xth$. Therefore, for the actual martingale $\ovS_t=\SG_t-t$ the corresponding  exponential martingale reduces trivially to the constant 1.  In contrast, the shifted $\symS-\delta t$ admits the value $c=\xth(\delta)$ for which $\ex{-c(\symS_t-\delta t)}$ is a martingale (solving the corresponding equation, e.g., for $\delta=1$, $c=0.714556...$ satisfies $c=\ln(1-c^2)$).
\vv
The $p$-moments are of the form $c_p^t$, with the base $c_p>1$ for $p>1$:
\[
c_p=c_{p,\xth}=\defif
{\Frac{(1+\xth)^{p}}{1+p\xth}}{p\xth>-1\,\,\text {for $M$}}
{\left(\Frac{(1-\xth^2)^p}{1-p^2\xth^2}\right)^{1/2}}{|p\xth|<1\,\,\text {for $\wt M$}}
\]
They are not uniformly bounded with the exception of $p=1$ which nevertheless entails an a.s.\ limit  as $t\to\infty$ by Doob's Convergence Theorem. However, these limits are 0 a.s.\ for $\xth\neq 0$. Indeed, e.g., denote $\SX _{\xth}=\lim_{t\to\infty} M_t(\xth)$. Also, by the SLLN $\SG_t/t\to 1$ a.s.. Denoting the event of convergence by $A_{\xth}$, we arrive at a contradiction:
\be\label{contr}
(1+\xth)e^{-\xth} =\lim_{t\to\infty} M_t^{1/t}(\xth)=1\quad\text{on}\quad
 A_{\xth}\cap \{\SX _{\xth}>0\},
\ee
so $\PR(\SX _{\xth}>0)=0$; similarly for  $\wt M_t$ with the SLLN  $\symS_t/t\to 0$. 
\vv
{\bf Quadratic variations}
\vv
The quadratic variation  $[M]$ contains nothing stochastic, just jumps.  That is, interpreting Notes \ref{phix} and \ref{ay}, putting $y=\varphi(x)$, 
with $x_t=\SG_t$ and $\varphi(x)=e^{-\xth x}$, $y=\varphi(x)$ and $a(t)=(1+\xth)^t$:
\be\label{dM}
\begin{array}{rcl}
de^{-\xth\SG_{t_-}} &=&-\xth ^{-\xth\SG_{t_-}}\, d\SG_t.\\
d \left[e^{-\xth\SG}\right]_t &=&\xth^2 e^{-2\xth\SG_{t_-}} \,d[\SG]_t; \\
dM_t &=&-\xth (1+\xth)^te^{-\xth\SG_{t_-}}\,d\SG_t+ \ln(1+\xth) (1+\xth)^t \,e^{-\xth\SG_{t_-}} dt,\\
     &=&-\xth M_{t_-}\,d\SG_t+ \ln(1+\xth) dM_{t_-} dt,\hspace{1.65in} (*)\\
d[M]_t&=&(1+\xth)^{2t}\, \xth^2 e^{-2\xth\SG_{t_-}} \,d[\SG]_t=\xth^2 M_{t_-}^2\,d[\SG]_t.\\
\end{array}
\ee
The compensator of $M^2$, a.k.a.\ the oblique bracket $\langle M\rangle$, stems from elementary stochastic calculus. Let $t>s$.
\[
\E[M_t^2|\ca G_s]= (1+\xth)^{2t} \E e^{-2\xth \SG_{t-s}} \, e^{-2\xth \SG_s}= \Frac{(1+\xth)^{2t}}{(1+2\xth)^{t-s}}.
\]
Hence
\[
\E[M_t^2-M_s^2|\ca G_s]= \left(\left[\Frac{(1+\xth)^{2}}{1+2\xth} \right]^t - \left[\Frac{(1+\xth)^{2}}{1+2\xth} \right]^s\right)
M_s^2 
\]
Denoting the base by $b_{\xth}=\Frac{(1+\xth)^{2}}{1+2\xth}$ for the sake of brevity, we obtain
\[
d\langle M,M\rangle_t=(\ln b_{\xth}) \,b_{\xth}^t \,M_{t_-}\,dt
\]
\vv
{\bf Induced fields}
\vv
Continuous linear operators w.r.t.\ to the variable $\xth$ preserve the martingale property subject to suitable integrability conditions. 
\vv
{\bf Example 1}. Partial derivatives $\frac{\partial^n}{\partial\xth^n}$.
\vv
Consider $M(\xth)=a(\xth)e^{-b\xth}$ with $a\in C^\infty$, and $b=\SG_t$ or $b=\symS_t$. By standard calculus
\[
M^{(n+1)} = e^{-b\xth}\sum_{j=0}^n (-1)^{n-j} {n\choose j}a^{(j)} b^{n-j}.
\]
The failing factorial simplifies the notation,
\[
(a)_j\df \Frac{\Gamma(a+1)}{\Gamma(a+1-j)}\qquad \text{(when $j\ge 1$, it's $\underbrace{a(a-1)\cdots (a-j+1)}_{j}$ )},
\]
entailing the binomial coefficient ${a\choose j} =\Frac{(a_j)}{j!}$. The case of $\SG_t$ is easy:
\[
M^{(n)}_t=M_t \Sum_{j=0}^n (-1)^j
{n\choose  j} (t)_{n-j} (1+\xth)^{-n+j} \SG_t^{n-j}.
\]
For $\symS_t$, the formula is similar but the coefficients are way more cumbersome. However, for
$a(\xth)=(1-\xth^2)^{t/2}$, using the binomial Taylor series,
\[
a^{(j)}(0)=\defif 0 {j=2m-1} {(-1)^m {t/2\choose m} (2m)!} {j=2m}.
\]
Plugging in also $\xth=0$ in the former expansion, 
we obtain polynomial martingales:
\[
\begin{array}{rcrcl}
P_n(\SG_t)&\df& (-1)^n M^{(n)}_t(0)&=&\Sum_{j=0}^n (-1)^{j} {n\choose  j} (t)_j  \SG_t^{n-j} ,\\
&&P_1(t)&=& \SG_t-t,\\
&&P_2(t)&=&\SG_t^2-2t\SG_t+(t)_2,\\
&&P_3(t)&=&\SG^3_t-3t\SG_t^2+3(t)_2 \SG_t-(t)_3,\\
\vspandexsmall
&&P_4(t)&=&\SG_t^4-4t\SG_t^3+6(t)_2\SG_t^2-4(t)_3\SG_t+(t)_4, \,...\\
\vspandexsmall
\wt P_n(\SG_t)&\df& (-1)^n\wt M^{(n)}_t(0)&=&\Sum_{2j\le n}(-1)^j {n\choose  2j} {t/2 \choose j} (2j)!\,\symS_t^{n-2j} ,\\
&&\wt P_1(t)&=& \symS_t,\\
&&\wt P_2(t)&=&\symS^2_t-t,\\
&&\wt P_3(t)&=&\symS^3_t-3t\,\symS_t,\\
&&\wt P_4(t)&=&\symS_t^4-4t\symS_t^2+6 t(t-2), \,...\\
\end{array}
\]
{\bf Example 2}. The Laplace transform of the function $\xth \mapsto M_t(\xth)$.
First,
\[
L_t(\theta)=\Frac{\Gamma(t+1)\,e^{{\SG}_t}}{({\SG}_t+\theta)^{t+1}}\quad\mbox{with}\quad \E L_t(\theta)=\frac{1}{\theta};
\]
then, through differentiation  w.r.t.\ $\theta$,
\[
(-1)^{k-1}\,L^{(k-1)}_t(\theta)=\Frac{\Gamma(t+k)\,e^{{\SG}_t}}{({\SG}_t+\theta)^{t+k}}\quad\mbox{with the mean}\quad\frac{(k-1)!}{\theta^k},\quad k\in\N.
\]
These martingales are not $p$-integrable for $p>1$ but  converge a.s.\ to 0. The argument is similar to one used in \refrm{contr}; i.e., by applying the exponent $1/(t+k)$, then SLLN and the Stirling's formula.
\vv
The case of $\symS_t$ seems very cumbersome since the Laplace transform of $(1-\xth^2)^{t/2}, \,|\xth|\le 1$ does not express in terms of standard functions. However, we may conceal the difficulty using Fubini's theorem:
\[
\wt E_\xi \df \ca E_{\xth} \wt S_t = \E' e^{-\xth (\symS_t-\symS'_t)},
\]
where the prime indicates an independent copy.  Thus, the transform $\wt L(\theta)=\Int_{-1}^1 \wt E_{\xi}e^{-\theta \xth}\,d\xth $ evaluated at $\theta=0$ yields a martingale
\[
\wt L_t=
2 \E'\left[\Frac{\sinh (\symS_t-\symS'_t)}{\symS_t-\symS'_t}\right]
\]

\appendix
\section{Appendix}
\setcounter{equation}{0}
\counterwithin*{equation}{subsection}
\renewcommand{\theequation}{\thesubsection.\arabic{equation}}
\subsection{Modular spaces}\label{modspaces}
We intend to describe the vector space of Poisson-integrable functions in terms of topological vector spaces but desirably  as F-spaces (cf.\ \cite{Rud}). In particular, the modular spaces (introduced in \cite{Mus}) become an essential tool in study of stochastic processes with independent increments, often quite elaborate in such generality, cf., e.g., \cite{KwaW}. This review befits diffusion-free or even trend-free L\'evy processes and thus benefits from their stationarity and simplified structure.
\vv
Let $(\BT,\ca T,\tau)$ be a $\sigma$-finite measure space entailing the F-space  $L^0(\BT)$ of measurable real functions, equipped with a metrizable topology of local convergence in measure. Consider the class $\cM$ of functions $\phi:[0,\infty)\to [0,\infty)$, $\phi\neq 0$, $\phi(0)=0$ that are 
\begin{enumerate}
\item
(a) nondecreasing and (b) continuous at 0;
\item
satisfying the inequality
\be\label{a+b}
\phi(\alpha u+(1-\alpha)v)\le\phi(u)+\phi(v),\quad \alpha\in [0,1].
\ee
\end{enumerate}
Then the integral
\be\label{mod}
\Phi(f)=\int _{\BT} \phi(|f|)\,d\tau
\ee
is a leading example of a modular on a metrizable {\bf\em modular space} (a.k.a.\ Musielak-Orlicz or generalized Orlicz space):
\be\label{Lmod}
L^\phi(\BT)=\set{f\in L^0(\BT): \Phi(f/c)<\infty \mbox{ for every $c>0$}},
\ee
which  becomes an F-space under the F-norm
\be\label{metric}
||f||_{\phi}=\inf\set{c>0: \Phi(f/c) \le c}.
\ee
 The completeness of the modular spaces \refrm{Lmod} follows from the continuity of the embedding $L^\phi(T)\hookrightarrow L^0(T)$ for any set $T\in\ca L$ of finite measure.
 \vv
Members of $\ca M$ are modulars, too ($\tau=\delta_1$, $L^\phi(\BT)=\R$).
\vv
A concave $\phi$ is subadditive; then  the modular $\Phi(f)$ itself is an F-norm.
The modification of a concave $\phi\mapsto \wt{\phi}(u):=\phi(u^2)$   preserves \refrm{a+b}  but may destroy subadditivity, so the resulting modular may fail the triangle inequality yet keep  a valid F-norm \refrm{metric}. We simply write $\wt{L}^{\phi}\df L^{\wt{\phi}}$.
\vv
Two functions are deemed equivalent, $\phi_1\sim \phi_2$,  if for $i,j\in\{1,2\}$
\be\label{simphi}
\exists~a,b:\, 0<b\le a~\forall~ u\ge 0 ~~ a\phi_i(bu)\le \phi_j(u).
\ee
Then $L^{\phi_1}=L^{\phi_2}$ as vector spaces. If $\phi_1$  is not a modular but $\phi_2$ is, then $L^{\phi_1}$ is metrizable with the help of the metric $\|f\|_{\phi_2}$ in spite of the functional $\|f\|_{\phi_1}$ failing to be  a metric. That is,  $\|f\|_{\phi_i}<\epsilon$ implies $\|f\|_{\phi_j}<\epsilon/b$
\vv
If one is a modular, we may call the other a modular (or a quasi-modular), too. Then the same symbol is kept (with some abuse of notation); e.g., $\phi_0(u)$ denotes  either of the following modulars:
\be\label{phi0}
u\wedge 1\quad\sim\quad \frac u {u+1}\quad\sim\quad 1-e^{-u}
\quad\sim\quad \arctan (u),\quad \text{etc.}
\ee
(or additional ones possibly failing monotonicity or \refrm{a+b}).
\vv
A one-sided L\'evy measure $\nu$ on $\R_+$  yields a perfect modular
\[
\phi^\nu(u)=\int_{\R_+} \left(1-e^{ux}\right)\,\nu(dx)\sim
\int_{\R_+} \phi_0(ux)\,\nu(dx).
\]
The need for extension of the concept of modulars is illustrated by troublesome but challenging functions:
\be\label{cosin}
u\,\mapsto \,\int_{\R_+} \left(1-\cos(ux)\right)\,\nu(dx),\quad
 \left|\int_{\R_+} \sin(ux)\,\nu(dx)\right|.
\ee
The latter function is dominated by the  modular
$\phi^\nu(u)=\int \R_+ (ux\wedge 1)\,\nu(dx)$ and the former one by $\wt{\phi}^\nu$. However, their modularity, even up the equivalence, is uncertain unless $\nu$ is special. E.g., for a $\alpha$ stable subordinator both functions are equivalent to $u^{\alpha}$. For a Gamma process we obtain
\be\label{phis}
\begin{array}{rcl}
\phi_1(x) &\df& \ln(1+x),\\
\phi_2(x)&\df& \frac 1 2 \ln(1+x^2)\sim \wt{\phi}_1.\\
\end{array}
\ee
Similarly, the function $\phi(u)=\left|\E\big(Ku: |Ku|\le 1\big)\right|$, although not being itself a modular in general but, when combined with the modular $\E\left(|xK|^2\wedge 1\right)$ may become one, subject to additional properties of the distribution $\kappa$, e.g., $u^\delta$ for $K$ with regularly varying tails with a parameter $\delta$, cf.\ \cite{Cli}.
\vv
For general processes with independent increments the integrability in terms of topological vector space becomes quite convoluted (cf.\ \cite[0.7-0.9]{KwaW}).
\eject
\subsection{Calculus of jumps}\label{jumps}
We consider  {\em pure jump processes} (PJP), mainly deterministic and only later suitably randomized as needed.
We write values of a function as $f(t)$ or $f_t$ as convenient.
\subsubsection{Basic modifications of functions}
Consider a countable set $\ca U=\{u_k:k\ge 0\}\subset [0,\infty)$ with  $u_0=0$. The ascending  rearrangement of a finite subset $\{u_1,..., u_n\}$ is denoted by $\wek u_n^*=(u_1^*,...,u_n^*)$ which defines the associated permutation $\sigma_1,\cdots,\sigma_n$ of indices: $u_{\sigma_k}=u_k^*$.
$\ca U$ entails the counting measure  $\nn =\sum_n \delta_{u_n}$:
\[
\begin{array}{l}
\nn  A=|A\cap U|=\Sum_n \I A(u_n),   \quad 
\nn (t)=\nn  [0,t] \text{\,\,(may be $\infty$ for all $t$)}, \\
\quad \nn  f\df\Int_0^\infty f(s)\,\nn (ds)=\Sum_n f(u_n),\quad \nn _t f=\nn  f\I{[0,t]}.\\
\end{array}
\]
We consider jumps $h_n$ with $h_0=0$ that may takes values in an F-space, e.g., in a modular space of random variables.
Under the mode of convergence, consider
\be\label{reqproof}
\left.\begin{array}{rl}
\vspandexsmall
\text{\bf the standing assumption:}&
\text{jumps are summable};\\
\text{\bf with a routine proof}: &\left\{
\begin{array}{l}
1.\, \text{Order linearly a finite $\ca U_0\subset \ca U$;}\\
2.\, \text{Use some discrete identity;}\\
3.\, \text{Finish by approximation.}\\
\end{array}\right.
\end{array}\right\}
\ee
Recall that a sequence $h_n$ is summable if the net of finite sums $\sum_{n\in N} h_n$ converges w.r.t.\ the natural partial order of inclusion. An F-space enjoys several equivalent conditions, studied as early as in 1930s by Orlicz (cf.\ \cite[Th.\ 1.3]{Szu} and references therein). In particular, the summability is equivalent to the permutation-invariant convergence and to the bounded convergence (i.e., with arbitrary bounded scalar multipliers). However, such conditions may fail to be equivalent with respect to a non-metric convergence such as the a.s.\ convergence (cf.\ \cite[Sect.5.3]{Szu} and references therein).
\vv
A right continuous {\em amassing\footnote{a.k.a.\ ``the cumulative function'' when $u_k$ are increasing} function} and integrals stem from {\em compounding}:
\be\label{xf}
\begin{array}{rcl}
xA&\df& \Sum_n h_n \I A(u_n),  \quad \text{hence}\quad x_t= \Sum_n h_n \Ib {u_n\le t} \quad\text{and}\\
xf&=&\Int_0^\infty f(s)\,x(ds)= \Sum_n h_n f(u_n);.
\end{array}
\ee
\begin{rem}\label{pure or not}\sf
In general, random jumps should be viewed as vectors in $L^0$ or another modular space of random variables. In other words, the associated measure is a vector measure rather than a mapping from the underlying probability space into the space of real (possibly signed) measures. Therefore, a purist may refuse to name the associated jump process ``pure''  with a possible exception of positive random jumps. \QED
\end{rem}
By the aforementioned Orlicz Summation Theorem, the integral is well defined for bounded real functions $f$ which entails the question of the maximal structure of integrable functions. Compounding yields a new PJP  $ y_t\df x_t f=xf\I{[0,t]}$ stemming from $x_t$ as $x_t$ stems from $\nn_t$).
\vv
We may display our quantities either as sums or as Lebesgue integrals  with respect to the discrete measure $x$, like in \refrm{xf}. For lcrl functions $f$ we may convey Riemann-Stieltjes integrals
$\int_0^t f_{s_-}\, dx_s$, then we may write concisely $dy_t=f_{t_-}\,dx_t$.
\vv
The old jumps:
\be\label{Delta}
\Delta x_t\df x_t-x_{t_-}=\defif {h_n}{t=u_n} 0{t\notin \ca U}.
\ee
produce new jumps by further compounding or compositions:
\be\label{xk} 
\begin{array}{rcl}
(x\circ \wek k)_t &\df& \Sum_n k_n h_n\Ib{u_n\le t}\\
           &=:& \left[x,y\right]_t,
               \,\,\text{where}\,\,  y_t=\Sum_n k_n \Ib{u_n\le t};\\
\left[x\right]^\varphi _t
&\df&  \Sum_n \varphi(h_n) \Ib{u_n\le t}=\lim_{\pi_t} \sum_k \varphi(x_{t_k}- x_{t_{k-1}}),
\end{array}
\ee
where $\pi_t=\{t_k\}$ are partially ordered partitions of $[0,t]$ with mesh $\to 0$ (having chosen a finite $\ca U_0\subset \ca U\cap [0,t]$, we may take $\pi_t \supset \ca U_0$). Clearly, the products $k_nh_n$ must be first well defined and must form a summable sequence.
\vv
The distinction between the operations is as valid  as the distinction between the verbs ``compound'' and ``compose''. Essentially, operations entail each other; the first operation by repetition yields $[x]^{(m)}$, where $\varphi(x)=|x|^m$; then entails the {\em variation} $[x]^\varphi$ at least when $\varphi$ is analytic. Conversely,  the {\em quadratic variation} $[x]^{(2)}$ with $\varphi(x)=x^2$ produces  $[x,y]=([x+y]^{(2)}-[x-y]^{(2)})/4$ by polarization. To distinguish $[x]$  from the {\em absolute variation}, 
we may denote the latter by $\|x\|\df[x]^{(1)}$ if necessary.
\vv
The change of the order of operations, i.e., $\varphi$ after the sum, may complicate a presentation of jumps and force the integral display. On the other hand, a classical integral representation may not exist; nevertheless, we can still use a symbolic differential notation.
For example, consider $\varphi(x)=x^2$: 
\be\label{squares}
\begin{array}{rcl c rcl}
\vspandex
[x]_t&=&
\Sum_n h_n^2 \Ib{u_n\le t} &\text{or}& 
d[x]_t&=&\Delta^2 x_t =\Delta x_t\,dx_t, \\
x^2_t&=& \Sum_m\Sum_n h_m h_n \Ib{U_n\vee U_m\le t}   &\text{or}&
dx^2_t&=& \Delta x_t^2 \,d\nn_t\\
 &&&& &=&   2 x_{t_-}\,dx_t+[x]_t
 \end{array}
 \ee
In general, one may expect a rather complicated presentation of jumps (e.g,, for the power $\varphi(x)=|x|_t^m$  cf.\ \cite[Sect.\ 9.2.3]{Szu}). 
 \begin{note} \label{phix}
 For the composition $y=\varphi\circ x$  with $\varphi\in C^1$:
\be\label{xdn}
\begin{array}{rcl}
\varphi(x_t)&=&\varphi(h_0)+\Int_0^t \Big(\varphi(x_s)-\varphi(x_{s_-})\Big)\,dx_s,
 \quad\text{or}\\
&&d\varphi(x)_t=\varphi'(x_{t_-})\, dx_t,\quad\text{also}\\
&&d[\varphi(x)]_t=|\varphi'(t)|^2 \,d[x]_t
\end{array}
\ee
\end{note}
\Proof We use of the telescoping sum
$\varphi_n=\varphi_0+\Sum _{k=1}^n (\varphi_k-\varphi_{k-1})$ for the integral formula, and then use the symbolic differential display.
\QED

The product $a_t y_t$ of a non-PJP function $a_t$ and a PJP $y_t$ is no longer PJP. 
\begin{note}\label{ay}
For a function $a\in C^1$,
\[
\begin{array}{rcll}
\vspandexsmall
d(a_t y_t) &=&a_t\,dy_t+a'_t\,y_{t_-} dt,\\
\vspandexsmall
d[ay]_t&=&a^2_{t}\, d[y]_t.\\
\end{array}
\]
\end{note}
\Proof Use the standard decomposition of increments and of increments of squares. In the latter case, only the sum with the first term doesn't vanish in the limit:
\vv
\hspace{30pt} $|a_ty_t-a_sy_s|^2=a_t^2(y_t-y_s)^2+y_s^2(a_t-a_s)^2 +2 a_ty_s (a_t-a_s)(y_t-y_s).$
\QED

\subsubsection{Inverse}
A rcll PJP $x_t$ is nondecreasing when $h_n\ge 0$ and thus
yields the rcll inverse
\be\label{app:inv}
\ix_v =x^{-1}(v) \df\inf\set{v: x_t> v}
\ee
Since the infimum is attained due to the right continuity, then
\[
\set{\ix_v>t}=\set{ x_t\le v}, \quad\mbox{or equivalently,}\quad \set{\ix_v\le t}=\set{x_t> v}.
\]
Hence the inverse is an involution on the rcll class: $(x^{-1})^{-1} =x$.
Similarly, the lcrl inverse
$ \ix(t_-)=\sup\set{t: x_t< v}$
is an involution on the lcrl class.

\begin{center}
\begin{tikzpicture}[scale=0.3]
    \draw[help lines, lightgray] (1,1) grid (11,13);
    \draw [<->,thick] (1,13) node (yaxis) [above] {$t$}
        |- (11,1) node (xaxis) [right] {$a$};
     \draw [very thick, blue](1,1)--(3,1);
     \fill[blue] (1,1) circle (6pt);
    \draw [very thick, blue](3,3)--(6,3);
    \draw [dashed, very thick, blue, red] (3,1)--(3,3);
    \draw (.9,3)--(1,3) node[left] {{\tiny $x_1$}};
    \draw (3,1)--(3,.9) node[below] {{\tiny $u_1$}};
    \draw (6,1)--(6,.9) node[below] {{\tiny $u_2$}};
    \fill[blue] (3,3) circle (6pt);
    \draw [very thick, blue](6,7)--(8,7);
     \draw [dashed, very thick, blue,  red] (6,3)--(6,7);
    \draw (.9,7)--(1,7) node[left] {{\tiny $x_2$}};
    \draw (8,1)--(8,.9) node[below] {{\tiny $u_3$}};
    \fill[blue] (6,7) circle (6pt);
    \draw [very thick, blue](8,9)--(9,9);
     \draw [dashed, very thick, blue,  red] (8,7)--(8,9);
    \draw (.9,9)--(1,9) node[left] {{\tiny $x_3$}};
    \draw (9,1)--(9,.9) node[below] {{\tiny $u_4$}};
    \fill[blue] (8,9) circle (6pt);
    \draw [very thick, blue](9,11)--(11,11);
    \draw [dashed, very thick, blue,  red](11,11)--(11,12.5);
     \draw [dashed, very thick, blue, red] (9,9)--(9,11);
    \draw (.9,11)--(1,11) node[left] {{\tiny $x_4$}};
    \fill[blue] (9,11) circle (6pt);
    %
    \draw[help lines,lightgray] (15,1) grid (27,11);
    \draw [<->,thick] (15,11) node (yaxis) [above] {$a$}
        |- (27,1) node (xaxis) [right] {$t$};
    \draw [very thick, blue, red] (15,3)--(17,3);
    \draw (14.9,3)--(15,3) node[left] {{\tiny $u_1$}};
    \draw (17,1)--(17,.9) node[below] {{\tiny $x_1$}};
    \fill[red] (15,3) circle (6pt);
    \fill[red] (17,6) circle (6pt);
    \draw [dashed, very thick, blue](17,3)--(17,6);
    \draw (14.9,6)--(15,6) node[left] {{\tiny $u_2$}};
    \draw [very thick, blue, red](17,6)--(21,6);
     \draw [dashed, very thick, blue] (21,6)--(21,8);
    \draw (14.9,8)--(15,8) node[left] {{\tiny $u_3$}};
    \draw (21,1)--(21,.9) node[below] {{\tiny $x_2$}};
    \fill[red] (21,8) circle (6pt);
    \draw [very thick, blue, red](21,8)--(23,8);
    \draw [dashed, very thick, blue](23,8)--(23,9);
    \draw [very thick, blue, red](23,9)--(25,9);
    \draw [dashed, very thick, blue](25,9)--(25,11);
    \fill[red] (23,9) circle (6pt);
    \draw (14.9,9)--(15,9) node[left] {{\tiny $u_4$}};
    \draw (23,1)--(23,.9) node[below] {{\tiny $x_3$}};
    \draw (24.9,1)--(25,.9) node[below] {{\tiny $x_4$}};
    \node[right] at (17,11.8) {\footnotesize rcll inverse};
\end{tikzpicture}

{$\begin{array}{c}
h_k>0\,\text{and}\,u_k\nearrow\\
x_k=h_1+\cdots+h_k,\,k=1,...,n
\end{array}$}
\end{center}

\begin{note}\label{app:inv:cont}
 If $\{h_n\}$ is summable and  $\ca U$ is dense then the inverse is continuous.
\end{note}
Indeed, w.l.o.g.\ we consider $[0,1]$. We truncate the sum for $x_t$ at the $N$\tss{th} term and sort $(u_1,\dots,u_N)\mapsto (u_1^*,\dots,u_n^*)$ in the ascending manner, adding $u_0^*=0$.  That is, with the corresponding permutation $\sigma$ of indices,
\[
{_{_N}}x_t=\Sum_{n=1}^N h_n \Ib{u_n\le t}=\Sum_{n=1}^N h_{\sigma(n)}\, \Ib{u^*_n\le t},
\]
Then, denoting $g_n=\Sum_{k=1}^n h_{\sigma(k)}$ and accumulating terms, for the inverse
\[
{_{_N}}\ix_t (v)=\sum_{n=1}^N (u_n^*-u_{n-1}^*)   \Ib{g_n\le v},
\]
the increments
$
{_{_N}}\ix_t-{_{_N}}\ix_{v_-}\le{\rm mesh}(u_n^*)\to 0
$. Hence $\Delta \ix_v=0$.
\QED

\addcontentsline{toc}{section}{References}
\vv\vv

\end{document}